\documentclass[final]{siamart251104}

\usepackage{lipsum}
\usepackage{amsfonts}
\usepackage{amssymb}
\usepackage{graphicx}
\usepackage{epstopdf}
\usepackage{algorithmic}
\ifpdf
  \DeclareGraphicsExtensions{.eps,.pdf,.png,.jpg}
\else
  \DeclareGraphicsExtensions{.eps}
\fi

\newsiamremark{remark}{Remark}
\newsiamremark{hypothesis}{Hypothesis}
\crefname{hypothesis}{Hypothesis}{Hypotheses}
\newsiamthm{claim}{Claim}

\headers{Quadratic Formula-based Nonlinear Approximation}{Z. He, C. Chen, M. H. Cho, J. Huang, and Y. Wu}

\title{Quadratic Formula-based Nonlinear Approximation\thanks{Submitted to the editors DATE.
}}

\author{Ziqin He\thanks{Department of Mathematics, University of North Carolina at Chapel Hill, Chapel Hill, NC 27599, USA (\email{ziqin@unc.edu}, \email{huang@email.unc.edu}).} 
\and Can Chen\thanks{School of Data Science and Society, Department of Mathematics, and Department of Biostatistics, University of North Carolina at Chapel Hill, Chapel Hill, NC 27599, USA (\email{canc@unc.edu}).}
\and Min Hyung Cho\thanks{Department of Mathematics and Statistics, University of Massachusetts Lowell, Lowell, MA 01854, USA (\email{minhyung\_cho@uml.edu}).}
\and Jingfang Huang\footnotemark[2] \thanks{Corresponding Author} \and Yichao Wu\thanks{Department of Mathematics, Statistics, and Computer Science, University of Illinois Chicago, Chicago, IL 60607, USA (\email{yichaowu@uic.edu}).} }

\usepackage{amsopn}

\makeatletter
\newcommand*{\addFileDependency}[1]{
  \typeout{(#1)}
  \@addtofilelist{#1}
  \IfFileExists{#1}{}{\typeout{No file #1.}}
}
\makeatother

\ifpdf
\hypersetup{
  pdftitle={Quadratic Formula-based Nonlinear Approximation},
  pdfauthor={Z. He, C. Chen, M. Cho, J. Huang, and Y. Wu}
}
\fi

\begin{document}

\maketitle

\begin{abstract}
  This paper presents a quadratic formula-based nonlinear representation for a given single-variable function $f(x)$, $-1 \leq x \leq 1$. First, we construct the explicit polynomial coefficient functions $a(x)$, $b(x)$, and $c(x)$ using a least-squares approach. Then, $f$ is reconstructed by solving the degree-2 polynomial equation $a(x) f^2 - b(x) f - c(x)=0$ for any $x \in [-1,1]$, where an index function is used to select the correct sign in the quadratic formula. The quadratic formula-based nonlinear approximation (degree-2 in $f$) outperforms classical orthogonal polynomial-based least-squares approximation (degree-0 in $f$) and rational approximation (degree-1 in $f$) for functions with sharp transitions or discontinuities.  As a potential application,  we apply the degree-2 representation to data denoising. Instead of relying on more complex ``edge-preserving'' metric-based optimization techniques, the smooth coefficient functions $a(x)$, $b(x)$, and $c(x)$ enable effective least-squares-based denoising on the low-dimensional manifold described by the algebraic variety $a(x) f^2 - b(x) f - c(x)=0$. 
 Denoising the index function, which determines the appropriate root to select, can be achieved using classical statistical or modern classification/clustering techniques. Numerical results and data denoising examples are provided to demonstrate the effectiveness of the degree-2 nonlinear approximation technique. The new nonlinear, quadratic formula-based representation also raises theoretical and numerical questions, including strategies for identifying numerically stable representations, developing optimal algorithms to construct the polynomial coefficient functions $a(x)$, $b(x)$, and $c(x)$, and achieving economical representation and denoising of the index function. Degree-2 representations of multivariable functions and higher degree-$d$ representations are currently under research, and results will be presented in the future.
\end{abstract}



\section{Introduction}
While linear combination-based function representations have been extensively studied 
in numerical analysis and optimization, nonlinear approximation theories and techniques 
remain relatively underexplored within the scientific computing community. 
Nevertheless, the effectiveness of nonlinear representations and operations has 
been well established. For example, convolutional neural network architectures employ 
nonlinear activation functions such as ReLU and composite functions across multiple 
layers as fundamental building blocks \cite{schmidhuber2015deep}. In Quadrature by Two Expansions (QB2X) \cite{ding2021quadrature,weed2023quadrature}, two-dimensional (2D) Laplace and Helmholtz single- and double-layer potentials are efficiently approximated using roots of a polynomial equation derived from the boundary geometry and the evaluation point. QB2X provides a more accurate representation of the layer potential within the leaf box of the fast multipole hierarchical tree structure compared to classical linear combination-based harmonic expansion schemes \cite{klockner2013quadrature,wala2019fast}.

In classical approximation theory, a given function $f(x)$ is expressed as a linear combination of basis functions $\phi_n(x)$, with the coefficients determined via the least-squares method. This approach is referred to in this paper as a degree-0 representation. The degree-0 representation can be extended by including an additional set of basis functions $\{ \phi_n(x) f(x) \}$, yielding the degree-1 representation of $f$. It can be shown that this degree-1 representation is equivalent to a rational function approximation. 

One of the main contributions of this paper is the extension of degree-0 and degree-1 representations into a quadratic formula-based degree-2 nonlinear approximation scheme. The degree-2 representation employs three sets of basis functions $\{ \phi_n(x) \}$, $\{ \phi_n(x)f(x) \}$, and $\{ \phi_n(x)f^2(x) \}$. The function $f(x)$ is implicitly defined by a quadratic equation with polynomial coefficients $a(x)$, $b(x)$, and $c(x)$, such that the function values lie on the manifold described by the following algebraic variety:
\begin{equation}
a(x) f^2(x) - b(x) f(x) - c(x) =0.
\end{equation}
To approximate the value of $f(x)$ at a sample point $x_0$, we first evaluate the coefficients $a(x)$, $b(x)$, and $c(x)$ at $x_0$, then apply the quadratic formula to compute the roots of the degree-2 representation at that point. The appropriate sign in the quadratic formula is chosen using an index function $\zeta(x)$, which takes the value of either $+1$ or $-1$. Therefore, the degree-2 representation enables accurate approximation of piecewise smooth functions with jumps, rapid transitions, and multi-valued behavior. Here, the separation of the smooth-coefficient manifold equation $a(x) f^2(x) - b(x) f(x) - c(x) =0$ and the index function $\zeta(x)$ is one of the interesting features of the degree-2 representation. This feature can be leveraged to denoise data with discontinuities. Instead of using an ``edge-preserving" metric \cite{jain2016survey,strong2003edge}, the smoothness of the coefficient functions in the manifold equation enables the use of traditional least-squares optimization techniques to recover the functions $a(x)$, $b(x)$, and $c(x)$. A separate denoising scheme is then applied to the index function $\zeta(x)$. We present numerical results for denoising a single-variable function perturbed by various types of noise satisfying specific statistical properties.

We organize this paper as follows. Sec.~\ref{sec:degree-0-1} presents the classical orthogonal polynomial approximation (degree-0 representation) and the rational approximation (degree-1 representation). Sec.~\ref{sec:deg2-definition} introduces the generalized degree-2 representation and discusses its properties. Sec.~\ref{sec:deg2-construction} explains how to numerically construct the degree-2 representation using either a greedy algorithm or a rank-revealing QR-based ranking method. Sec.~\ref{sec:numcompare} compares the performance of various representations for different types of functions.  Sec.~\ref{sec:denoise} presents an application of the degree-2 representation to data denoising. Finally, Sec.~\ref{sec:conclusions} summarizes our results and outlines several ongoing projects aimed at further exploring this polynomial equation root- and manifold-based nonlinear approximation technique.

\section{Classical Degree-0 and Degree-1 Representations}
\label{sec:degree-0-1}
In classical interpolation and approximation theory, $f(x)$ is approximated by a 
linear combination of basis functions $\{\phi_n(x)\}_{n=0}^N$, where the coefficients
$\{ c_n \}$ are computed using least-squares optimization
\begin{equation}
\label{eq:deg0}
\min_{\{c_n\}} \left\| f(x) - \sum_{n=0}^N c_n \,\phi_n(x) \right\|_2 .
\end{equation}
The selection of basis functions has been extensively studied. In particular, when orthogonal bases, such as the Fourier series, 
wavelet families, and Chebyshev, Legendre, Laguerre, and Hermite polynomials are employed, the coefficients are computed as 
$$c_n = \frac{\langle f(x), \phi_n(x) \rangle}{\langle \phi_n(x), \phi_n(x) \rangle},$$
where  $\langle f(x), \phi_n(x) \rangle= \int w(x)f(x)\phi_n(x) dx$ is the inner product of $f(x)$ with the corresponding basis function under an appropriate weight function $w(x)$.
In this paper, we use the truncated non-orthogonal basis set $\{ x^n \}_{n=0}^{N}$ to illustrate the main ideas for simplicity. However, the normalized Legendre polynomial basis $\{ L_n(x)\}_{n=0}^N$, $x \in [-1,1]$ is used for most numerical experiments. It is well-known that the normalized Legendre polynomials are orthonormal to each 
other under the inner product with weight function $w(x)=1$, i.e., 
$\langle L_m(x), L_n(x) \rangle = \delta_{m,n}$. 
The truncated expansion provides the best mean-square ($L_2$) approximation of $f$ on $[-1,1]$, 
with expansion coefficients given by $c_n=\langle f,L_n\rangle$. We formally define the degree-0 representation as follows: 
\begin{definition}[Degree-0 representation]
We call the classical representation using the linear combination of basis functions 
$$\,f(x)\approx\sum_{n=0}^N c_n\,\phi_n(x)\,$$
a \emph{degree-0} representation because no (power of) $f$ appears in the basis.
The basis functions are explicitly given and only depend on $x$. When an orthonormal 
basis set with inner product weight function $w(x)$ is used, the expansion coefficients 
are $c_n=\int_a^b w(x)\,f(x)\,\phi_n(x)\,dx$.
\end{definition}

Motivated by one of the main ideas from the randomized linear algebra technique, where a matrix $\mathbf A$ is applied to several random vectors $\vec{v}_i$ and the results are used to extract compressible features of $\mathbf A$ \cite{martinsson2011randomized,martinsson2020randomized}, we can also include a function $f(x)$ multiplied by a simple polynomial function in $x$, in addition to the classical polynomial basis in Eq.~(\ref{eq:deg0}), as part of the least-squares optimization formulation. In other words, two sets of basis functions $\{1,x,\dots,x^{N_0} \}$ and $\{x\cdot f(x),x^2\cdot f(x),\dots,x^{N_1}\cdot f(x)\}$ are employed in the least-squares problem
\begin{equation}
\label{eq:deg1}
\min_{\{c_n,b_n\}} \left\| f(x) - \Big(\sum_{n=0}^{N_0} c_n x^n \;+\; \sum_{n=1}^{N_1} b_n\, x^n f(x)\Big) \right\|_2.
\end{equation}
Here we use the non-orthogonal monomials $\phi_n(x)=x^n$ for notational simplicity. Generalization to an orthonormal basis is straightforward and is often preferred in our numerical experiments because of its improved stability properties. Unlike randomized linear algebra algorithms, the basis functions $\{1,x,\dots,x^N\}$ are ordered 
by importance, which guides the truncation of the expansion. This formulation is referred to as the degree-1 representation:

\begin{definition}[Degree-1 representation]
We define
$$
f(x) \approx \Big(\sum_{n=0}^{N_0} c_n x^n \;+\; \sum_{n=1}^{N_1} b_n\, x^n f(x)\Big)
$$
as a \emph{degree-1} representation, since the degree-1 term $f(x)$ ($f$ power $1$)  explicitly 
appears in the second summation of the least-squares approximation in Eq.~\eqref{eq:deg1}.
\end{definition}

Solving this linear equation for $f(x)$ yields a rational approximation 
\cite{herremans2023resolution,trefethen2024polynomial,trefethen2025rational} 
\begin{equation}
\label{eq:rational}
f(x) \approx \frac{\sum_{n=0}^{N_0} c_n x^n}{1 - \sum_{n=1}^{N_1} b_n x^n} = \frac{c(x)}{b(x)},
\end{equation}
where $c(x)=\sum_{n=0}^{N_0} c_n x^n$ and $b(x)=1 - \sum_{n=1}^{N_1} b_n x^n$ are 
polynomial functions. 
Note that classical rational approximation solves a nonlinear optimization problem such as 
$\min_{P, Q}\,\|f - P/Q\|_{L_2}$ or $\min_{P, Q}\,\|f - P/Q\|_{L_{\infty}}$, whereas the degree-1 
formulation solves a linear optimization problem $\min_{b(x),c(x)} \| b(x)f(x) - c(x) \|_{L_2}$. 
To remove the trivial scaling ambiguity $c(x)/b(x)=(\alpha c(x))/(\alpha b(x))$, the current degree-1 representation enforces $b(0)=1$. The AAA algorithm \cite{nakatsukasa2018aaa} is one of the recent developments for constructing rational approximations.  It computes an approximate solution to the nonlinear optimization problem by iteratively solving a sequence of linear least-squares optimization problems with adaptively updated barycentric weights. 

The asymptotic convergence properties of degree-0 and degree-1 approximations are summarized for three classes of functions: (i)$f$ is analytic, (ii) $f$ has branch point singularities, and (iii) the more general case where $f$ or its derivatives exhibit discontinuities. We refer interested readers to \cite{trefethen2024polynomial,trefethen2025rational} 
for detailed analysis and comparisons of polynomial and rational approximations.
We first consider the degree-0 representation, which is equivalent to the ``polynomial
approximation" in \cite{trefethen2024polynomial}. When $f$ is analytic on  $[-1,1]$ and can be analytically extended to a neighborhood of $[-1,1]$, the $N$-term degree-0 representation 
converges exponentially at a rate of $R^{-N}$ for some $R>1$, the {\it analyticity radius of $f$} \cite{trefethen2024polynomial}.  When $f$ has a branch-cut singularity, e.g., $f(x)=|x-x_0|^\alpha$ with $x_0\in(-1,1)$,
the degree-0 representation only converges algebraically (polynomial error decay rate). In particular,
if $f(x)$ is discontinuous at a point $x_0 \in [-1,1]$, the Gibbs phenomenon occurs,
with approximation error overshoots and oscillations near $x_0$  \cite{gottlieb1997gibbs,hewitt1979gibbs}. 

When the condition $b(x)=1$ is enforced (or $N_1=0$), the degree-1 approximation reduces to the degree-0 
approximation. Therefore, the degree-1 representation can be considered as a generalization of
the degree-0 representation. When $f(x)$ is analytic on $[-1,1]$, it is straightforward
to conclude that the error of the ``optimal'' degree-1 approximation also converges at least exponentially 
at a rate of $R^{-N}$, where  $N=N_0+N_1$ and $R>1$.  An interesting case occurs when the approximation is required on a non-convex domain and $f(x)$ has a singularity in an ``inlet''
(see examples in \cite{trefethen2024polynomial,trefethen2025rational}). In this scenario, the analyticity radius of $f(x)$ is slightly greater than $1$ but very close to it.
This means that although convergence is exponential, the number of additional terms required to gain one significant digit is numerically impractical for the degree-0 representation. However, when the degree-1 representation is 
used and an optimal $b(x)$ is chosen, the analyticity radius of the function $b(x) f(x)$ increases, leading to much faster convergence of the error $\| b(x) f(x)- c(x) \|$. A larger $R$ value implies that the rational approximation will asymptotically outperform the degree-0 polynomial approximation. As shown in \cite{newman1964rational,trefethen2024polynomial,trefethen2025rational}, the degree-1 rational approximation can achieve
root-exponential convergence when $f(x)$ has branch point singularities, significantly outperforming the degree-0 polynomial approximation, which
converges only algebraically. When the maximum error is measured only over the region $[-1,-l] \cup [l,1] $ for any $1>l>0$, Zolotarev rational functions can be constructed analytically to approximate the sign function, achieving exponential convergence \cite{nakatsukasa2016computing,petrushev2011rational}.
Unfortunately, similar to the degree-0 representation and some classical rational approximations, if the function $f(x)$ is discontinuous at a point $x_0 \in [-1,1]$, the Gibbs phenomenon may still occur with the degree-1 approximation. The Gibbs phenomenon can be mitigated using a modified Pad\'{e}-Chebyshev approximation \cite{tampos2012accurate}.

\section{Degree-2 (Quadratic) Representation}
One of the main contributions of this paper is the generalization of the degree-0 and degree-1 
representations. A new set of basis functions $\{x^n f^2(x)\}_{n=1}^{N_2}$ is added to
existing sets $\{x^n\}_{n=0}^{N_0}$ and $\{x^n f(x)\}_{n=0}^{N_1}$, and $f^2(x)$ is approximated by solving the following least-squares problem:
\begin{equation}
\label{eq:deg2}
\min_{\{c_n,b_n,a_n\}}
\left\|\,f^2(x) -
\Big(\sum_{n=0}^{N_0} c_n x^n
+\sum_{n=0}^{N_1} b_n x^n f(x)
+\sum_{n=1}^{N_2} a_n x^n f^2(x)\Big)\right\|_{2}.
\end{equation}
Here, instead of directly approximating $f(x)$,
$$
f(x) \approx \sum_{n=0}^{N_0} c_n x^n \;+\; \sum_{n=0}^{N_1} b_n x^n f(x)
\;+\; \sum_{n=1}^{N_2} a_n x^n f^2(x),
$$
$f^2(x)$ is approximated to avoid potential issues related to existence, uniqueness, and other regularity concerns - for example, when $f$ is the sign function, it satisfies the smooth-coefficient manifold equation $f^2(x)\equiv 1$ without a linear $f$ term. 

\subsection{Definition and Properties of Degree-2 Representation}
\label{sec:deg2-definition}
Once the polynomial coefficients $a_n$, $b_n$, and $c_n$ that minimize Eq.~(\ref{eq:deg2}) are determined, a simple algebraic manipulation yields a manifold equation with smooth coefficients:
$$\left(1-\sum_{n=1}^{N_2} a_n x^n\right)f^2(x) - \left(\sum_{n=0}^{N_1} b_n x^n\right) f(x)-\left(\sum_{n=0}^{N_0} c_n x^n\right)\approx 0, $$
which serves as an implicit representation of $f(x)$. We define the {\it degree-2 representation} as follows.
\begin{definition}[Degree-2 representation]
\label{def:deg2}
The degree-2 implicit representation of a function $f(x)$ consists of two components: 
(a) {\it smooth} polynomial coefficient functions 
$$a(x)=1- \sum_{n=1}^{N_2} a_n x^n,~ 
b(x)=\sum_{n=0}^{N_1} b_n  x^n, \mbox{ and }~ c(x)=\sum_{n=0}^{N_0} c_n  x^n $$ 
describing the smooth-coefficient manifold $ a(x) f^2(x) - b(x)f(x) - c(x) =0$ (an algebraic variety) 
on which $f(x)$ lies; and (b) a single-bit index function $\zeta(x)$, which is 
either $-1$ or $1$, indicating which root should be selected. 
\end{definition}
For a single-valued function $f(x)$, once the coefficient functions $a(x)$, $b(x)$, $c(x)$, and the 
index function $\zeta(x)$ are determined, 
applying quadratic formula to the manifold equation $a(x) f^2(x) - b(x) f(x) -c(x) \approx 0$ gives: 
\begin{equation}
    \label{eq:quadratic}
 f(x) \approx \frac{b(x) + \zeta(x) \sqrt{b^2(x)+4 a(x) \cdot c(x)}}{2 a(x)}.
\end{equation}

The degree-2 representation generalizes the degree-0 and degree-1 
representations from the optimization perspective. By separating the manifold equation $a(x) f^2(x) - b(x) f(x) -c(x) = 0$ from the index function $\zeta(x)$, the degree-2 representation enables more effective algebraic manipulation of single-valued functions 
with discontinuities or data lying on low-dimensional manifolds. Sec.~\ref{sec:denoise} illustrates how this separation facilitates more effective data denoising. However, unlike the degree-0 and degree-1 representations, the degree-2 representation for single-valued functions requires additional storage and specialized tools to process the index function $\zeta(x)$, such as statistical or clustering/categorization techniques. 
An interesting observation is that all classical Walsh wavelet functions 
\cite{beauchamp1976walsh,walsh1923closed} lie on the same
manifold defined by the algebraic variety $f^2(x)-1=0$.  They can only be distinguished by
analyzing the index function $\zeta(x)$ for each basis, which is recursively defined using a hierarchical 
tree structure across different Walsh basis levels. This tree-based index function definition may offer more storage-efficient and effective operations compared to a point-wise bit-index function definition.

We now briefly analyze the convergence properties of the degree-2 representation. Since the degree-2
representation generalizes the degree-0 and degree-1 approximations, for analytic
functions, the ``best'' approximation error decays at least exponentially at a rate of $e^{-cN}$ for
some constant $c>0$, where $N=N_0+N_1+N_2$. For functions with branch-point singularities, one can expect at least
root-exponential convergence. The most notable 
feature of the degree-2 representation is its ability to effectively describe data with discontinuities or discontinuous derivatives. As a simple example, consider 
the ReLU activation function, $f(x)=\max(0,x)$,  widely used in machine learning. In this case, the degree-2 representation  is given by $a(x)=1$, $b(x)=x$, and $c(x)=0$, 
corresponding to the manifold $f(x) (f(x)-x) = f^2(x) - x \cdot f(x) = 0$. The index function
is defined as $\zeta(x) = -1$ if $x \leq 0$ and $\zeta(x) = 1$ otherwise.

For the more general case, let $f$ be defined as $f_-$ on $[-1, x_0]$ and $f_+$ on $[x_0, 1]$, where $-1< x_0 < 1$, and $f_-$ and $f_+$ are respectively analytic with $f_-(x_0) \neq f_+(x_0)$. Both $f_{\pm}$ can be independently approximated on each interval by the degree-0 representation with polynomial degree $N$,  yielding $p_{\pm}$ satisfying
\[
\|f_- - p_-\|_{L_2([-1,x_0])}\le C_-\,\rho_-^{-N},\qquad
\|f_+ - p_+\|_{L_2([x_0,1])}\le C_+\,\rho_+^{-N},
\]
where $C_{\pm}$ are positive constants and $\rho_{\pm}$ denote analyticity radii. Now define
\[
b:=p_-+p_+ ~\mbox{ and } ~c:=-\,p_-\,p_+.
\]
Then $b$ is a degree-$N$ polynomial, and the degree of $c$ can be reduced to approximately $N$ after truncation with a prescribed error tolerance, owing to the exponential error decay of both $p_-$ and $p_+$.  Finally, using  $b$ and $c$, the algebraic manifold equation of the degree-2 (quadratic) representation 
satisfies the exact algebraic identity
\[
f^2 - b f - c \;=\; (f-p_-)\,(f-p_+) = 0 \quad\text{on }[-1,1].
\]
Assuming the index function $\zeta(x)$ is correctly defined, the degree-2 representation is given by
\[
\hat f(x) :=
\begin{cases}
p_-(x), & x<x_0,\\[1ex]
p_+(x), & x>x_0,
\end{cases}
\]
which  converges exponentially over the entire domain $[-1,1]$,
\[
\|f-\hat f\|_{L_2([-1,1])} 
\;\le\; C\big(\rho_-^{-N}+\rho_+^{-N}\big). 
\]

Note that one can also use the degree-1 approximation instead of the degree-0 representation on
each subinterval. A larger analyticity radius would further improve the convergence of the resulting degree-2 representation. The optimized 
degree-2 representation is expected to outperform the degree-0 and degree-1 representations in terms of the number of basis functions required to achieve the same accuracy. 

\begin{figure}[htbp]
\centering
\begin{tabular}{cc}
\includegraphics[scale=0.5]{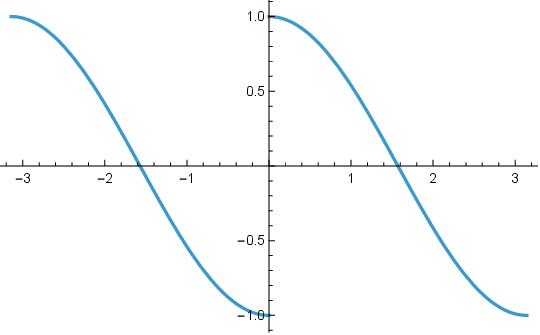}
&
\includegraphics[scale=0.5]{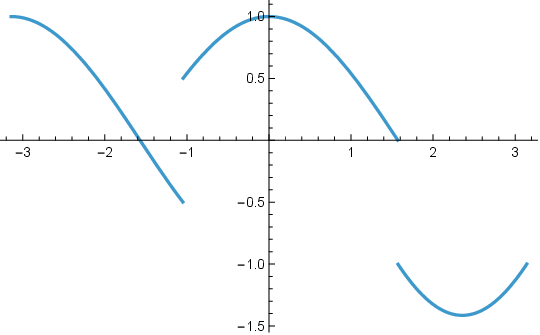} \\
\includegraphics[scale=0.5]{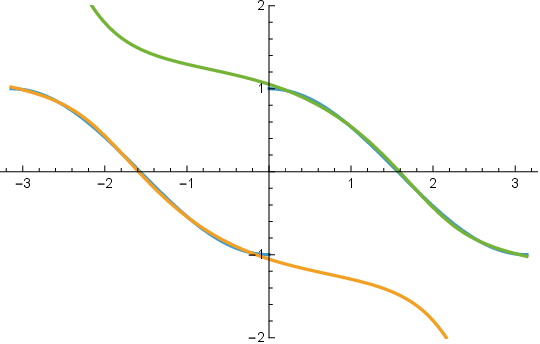} 
&
\includegraphics[scale=0.5]{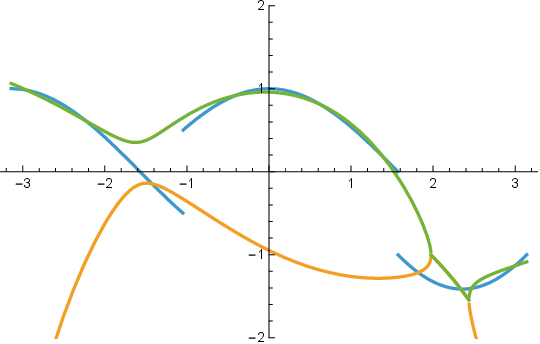} 
\end{tabular}
\caption{$f(x)$: \textcolor{blue}{Blue}. Root 1: \textcolor{orange}{Orange}, Root 2: \textcolor{green}{Green}. Top Left: $f(x)$ 
with one discontinuity. Bottom Left: well-separated branches of the manifold. 
Top Right: $f(x)$ with two jumps. Bottom Right: more complex branches
and index function for degree-2 representation.}
\label{fig:figJumps}
\end{figure}

As a novel approximation strategy, the degree-2 representation exhibits many interesting features, and further investigation is needed to establish its theoretical foundations. For example, Fig.~\ref{fig:figJumps} illustrates how the manifold properties change when a degree-2 representation is applied to $f(x)$ with one and two jump discontinuities. 
In the top-left panel of Fig.~\ref{fig:figJumps},  
$f(x)=\cos(x) (2 H(x) - 1)$, $x\in[-\pi, \pi]$, where $H(x)$ denotes the standard Heaviside function, is plotted in \textcolor{blue}{blue} color. The 
degree-2 approximation is computed as the optimal solution of Eq.~(\ref{eq:deg2})
with $N_2=0$ and $N_0=N_1=4$. In the bottom-left panel of Fig.~\ref{fig:figJumps}, the two curves (\textcolor{green}{green} and \textcolor{orange}{orange}) representing the two roots (manifold branches) are well-separated, each capturing one of the two smooth segments of the function. In the top-right panel of Fig.~\ref{fig:figJumps}, 
$f(x)=\left(2H(x+\frac{\pi }{3})-1\right) \cos (x)
   -H \left(x-\frac{\pi }{2}\right)  \sin (x),$ which has two discontinuities, is plotted. The behavior of the degree-2 representation (with $N_2=0$ and $N_0=N_1=4$) becomes more complex in this case. In the bottom-right panel of Fig.~\ref{fig:figJumps}, the two branches of the manifold attempt to jointly reduce the approximation error. This suggests that the index function may change value at points where $f(x)$ is smooth (e.g., at around $x\approx-1.7$), in addition to the discontinuity at $x=-1$. These joint efforts of the two manifold branches to reduce the approximation error may enable effective use of the degree-2 representation for functions with multiple jump points, where one branch can use the interval from another branch as a buffer region for a smooth transition to the next approximation interval.

\subsection{Numerical Construction of Degree-2 Representation}
\label{sec:deg2-construction}
As noted in Sec.~\ref{sec:degree-0-1}, we use normalized Legendre 
polynomials $\{ L_n(x) \}$ to construct the smooth-coefficient manifold equation $a(x) f^2(x) - b(x) f(x) - c(x)=0$ for improved numerical stability. The main task is to compute $a(x)=L_0(x)- \sum_{n=1}^{N_2} a_n L_n(x) $, 
$b(x)=\sum_{n=0}^{N_1} b_n L_n(x)$, and $c(x)=\sum_{n=0}^{N_0} c_n L_n(x)$ using the least-squares optimization in Eq.~(\ref{eq:deg2}).
We first sample the function $f(x)$ at a sufficiently large set of Gauss quadrature points $\{ x_i \}$, $i=1,\cdots, M$, which are the zeros of the Legendre polynomial $L_M(x)$. For numerical examples, we use $M=1000$. We define the degree-2 dictionary as the set of column basis vectors forming the matrix $\mathbf{V}$ 
\[
\mathbf{V} \;=\; 
\big[L_0,\dots,L_{N_0}\big]
\;\cup\;
\big[f\odot L_0,\dots,f\odot L_{N_1}\big]
\;\cup\;
\big[f^2\odot L_1,\dots,f^2\odot L_{N_2}\big],
\]
where $L_k$, $f$, and $f^2$ denote column vectors containing the Legendre polynomial $L_k$, $f$ and $f^2$ evaluated at $\{x_i\}_{i=1}^M$, respectively. The operator $\odot$ denotes the element-wise product of two vectors. The constant term $f^2 \odot L_0$ is omitted to avoid duplication.

Note that the inner product ($w(x)=1$) for continuous basis functions can be accurately 
computed using the Gauss quadrature rule:
\[
\langle u,v\rangle_w:= \int_{-1}^1 u(x)\,v(x)\,w(x)\,dx \;\approx\; \sum_{i=1}^M w_i\,u(x_i)\,v(x_i),
\]
where $w_i$ denotes the Gauss quadrature weight.
Define $\mathbf{W}:=\mathrm{diag}(w_1,\dots,w_M)$. The discretized least-squares problem is then formulated as
\[
\min_{c\in\mathbb{R}^K}\;\big\| \mathbf{W}^{1/2}(\mathbf{V} \eta - y)\big\|_2,
\]
where $\eta$ contains the coefficients of the functions $[a(x), b(x), c(x)]$, $y$ is the vector $ f^2 \odot L_0$, and $K=N_0+N_1+N_2+2$ denotes the total number of unknown coefficients.

Unsurprisingly, the choice of $N_0$, $N_1$, and $N_2$ significantly impacts the effectiveness of the degree-2 representation. A naive choice is to set $N=N_0=N_1=N_2$. However, as shown in Sec.~\ref{sec:numcompare}, this uniform selection strategy often leads to a larger number of basis terms in the representation for a given error tolerance. In the following, we present two approaches to improve the effectiveness of the degree-2 representation. 

\vspace{0.1in}
{\noindent \bf Greedy Algorithm for Basis Function Selection.} 
We first introduce a greedy algorithm to select additional basis functions and adaptively increase  
$N_0$, $N_1$, and $N_2$. Consider the following three candidate streams:
\[
\mathsf{S}_1 = [L_0,\dots,L_{n}],\quad
\mathsf{S}_2 = \mathsf{S}_1\odot f,\quad
\mathsf{S}_3 = [L_1,\dots,L_{n}]\odot f^2
\]
for a sufficiently large $n$, which can be adjusted if needed. Assume that 
at step $k$, the first $N_0$, $N_1$, and $N_2$ basis functions have already been selected from
$\mathsf{S}_1$, $\mathsf{S}_2$, and $\mathsf{S}_3$, respectively, and collected in the matrix $\mathbf{\Phi}_{k}$. Next, let $d$ be a set of 1, 3, or 5 basis vectors drawn from unused (and ordered) vectors in each candidate stream, and normalize them as
\[
  \widehat d \;=\; \frac{d}{\|d\|_{W}},
  \]
where $\|d\|_W := \sqrt{d^\top W d}$. Then, $\widehat{d}$ is tentatively added to $\mathbf{\Phi}_{k}$ to form a temporary basis set 
$ \mathbf{\widehat\Phi}_{k}=[\mathbf{\Phi}_{k}\;\,\widehat d]$. For each $\widehat{d}$,
solve the least-squares problem  
  \(
  \min_c \|W^{1/2}( \mathbf{\widehat\Phi}_{k} \eta - y)\|_2
  \)
and compute the residue 
$r=\|W^{1/2}(\mathbf{\widehat\Phi}_{k}\eta - y)\|_2$. Let $\widehat d^\star$ be the candidate that produces the smallest residue. Then, form the new basis set $\mathbf{\Phi}_{k+1}=[ \mathbf{\Phi}_{k}\;\,\widehat d^\star]$ for step $k+1$. If the residues $r$ are equal, randomly select $\widehat{d}$ from one of the three candidate streams. Finally, repeat this process until the residue $r$ meets the prescribed accuracy requirement. A QR decomposition can be used to recycle previous least-squares results and accelerate the solution process.

\vspace{0.1in}
{\noindent \bf Rank-revealing Basis Selection.} 
It is also possible to preprocess and select the most important basis vectors/functions from the candidates $[\mathsf{S}_1, \mathsf{S}_2, \mathsf{S}_3]$
using a rank-revealing QR factorization \cite{gu1996efficient}. Note that since an $m$-point Gauss quadrature rule can integrate polynomials of degrees up to $2m-1$ exactly, when the number of sample points $m$ exceeds $n$,
the vectors $q$ in $\textsf{S}_1$, representing sampled normalized Legendre polynomial values at Gauss quadrature nodes, form a set of orthonormal basis with respect to the $W$ norm, i.e., 
$q^\top W q = I$. 

Define $\mathbf{A}=\mathbf{W}^{1/2} [\mathsf{S}_1, \mathsf{S}_2, \mathsf{S}_3] $ and apply the rank-revealing QR algorithm:  
$\mathbf{A}\,\mathbf{\Pi} = \mathbf{Q_D R_D}$. For a given error tolerance, let $\mathbf{\Phi}$ denote the truncated $\mathbf{Q_D}$ and solve the
least-squares problem:
\[
  \min_{c}\; \big\|\mathbf{W}^{1/2}(\mathbf{\Phi} \eta - y)\big\|_2.
  \]
Finally, $a(x)$, $b(x)$, $c(x)$ can be recovered from $\eta$, $\mathbf{R_D}$, and the permutation matrix $\mathbf{\Pi}$.

\vspace{0.1in}
{\noindent \bf Comments.} Both the Greedy algorithm and the rank-revealing QR 
basis ranking method are not yet optimal for finding a smooth-coefficient manifold described by 
$a(x) f^2(x) - b(x) f(x) - c(x)=0$,
on which the data/function $f(x)$ lies. The existence of the relation 
$a(x) f^2(x) - b(x) f(x) - c(x)=0$ implies that the candidate set 
$[\mathsf{S}_1, \mathsf{S}_2, \mathsf{S}_3]$ should be treated as a (linearly dependent) 
frame set for constructing an improved compressed representation of $f(x)$. This relation equation 
is clearly not unique. It is well-known that the degree-0 representation is numerically backward stable; thus, 
the basis functions in $S_1$ are preferred for better stability properties. Proper ordering of the candidate vectors
in $[\mathsf{S}_1, \mathsf{S}_2, \mathsf{S}_3]$ ensures that $a(x)$ stays away from $0$ and that the quadratic formula 
in Eq.~(\ref{eq:quadratic}) remains numerically stable. We are currently generalizing the rank-revealing approach and exploring appropriate weights for each vector in $[\mathsf{S}_1, \mathsf{S}_2, \mathsf{S}_3]$ to obtain the most stable and efficient degree-2 representation.

\section{Numerical Comparison of Degrees-0, Degree-1, and Degree-2 Approximations}
\label{sec:numcompare}
The degree-2 representation generalizes the degree-0 and degree-1 representations, enabling the study of multi-valued and discontinuous functions on smooth-coefficient manifolds defined
by the roots of a polynomial equation (an algebraic variety). In this section, we numerically 
examine the performance of the degree-2 representation for: (1) a discontinuous function, 
(2) a highly oscillatory function, and (3) a function with a rapid transition region.

\vspace{0.1in}
{\noindent \bf Approximation of a Discontinuous Function.}
When \(f\) has a jump at a point in \([-1,1]\), classical polynomial (degree-0) expansions exhibit 
the Gibbs phenomenon. The coefficients of the Legendre polynomial expansion decay only algebraically at the rate of 
$|c_n|=\mathcal{O}(n^{-1})$. The degree-1 representation can significantly improve convergence. The degree-2
representation naturally places the two smooth pieces on separate curves describing a manifold defined by a polynomial equation
with smooth coefficients, and its manifold equation approximation can achieve global exponential convergence. 

\begin{figure}[htbp]
\centering
\includegraphics[width=8cm, height = 5cm]{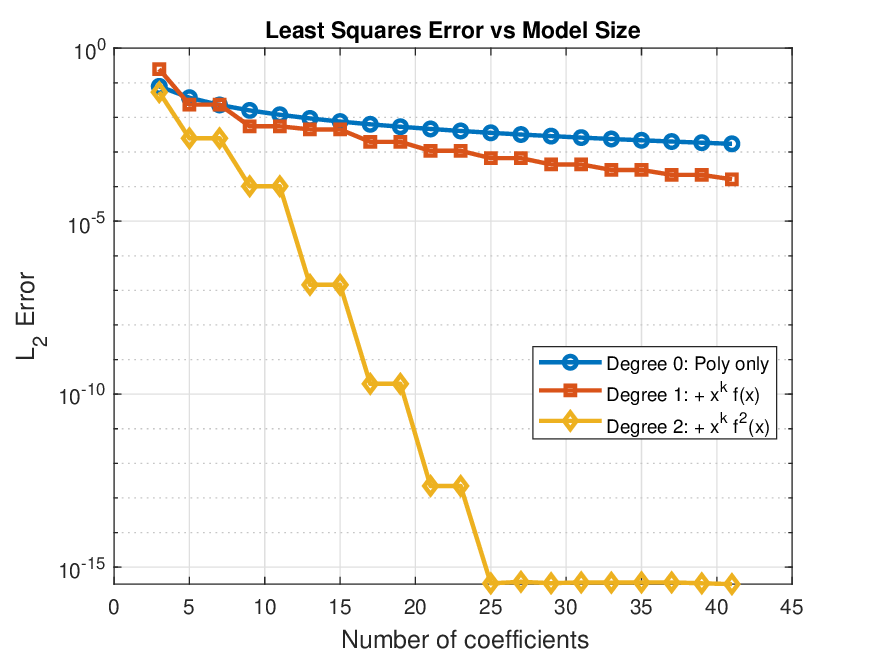}
\caption{Convergence of degree-0, degree-1, and degree-2 representations.
\label{fig:discontinuous} }
\end{figure}

In Fig.~\ref{fig:discontinuous}, we compare the convergence of 
these representations for the function $f(x)=\sin(x) (2 H(x) -1)$, 
$x \in [-1,1]$, where the Heaviside function $H(x)$ is defined as $0$ when $x < 0$ and $1$ when  $x \geq 0$. We use non-adaptive degree-1 ($N_0=N_1$) 
and degree-2 ($N_0=N_1=N_2$) representations. The $x$-axis shows 
the total number of coefficients in each representation (equal to the number of basis functions used),
and the $y$-axis shows the $L_2$ approximation error. As expected, the degree-0 
and degree-1 approximations converge slowly in this case. In contrast, the degree-2 
representation converges exponentially to machine precision.

In Sec.~\ref{sec:deg2-definition}, another example (Fig.~\ref{fig:figJumps}) illustrates how the degree-2 representation approximates a function with discontinuities. It shows how different branches of the manifold work together to reduce approximation error and how the index function may change value at points where the function is otherwise smooth. 
This latter feature is likely necessary to achieve optimal accuracy when the number of basis functions is insufficient. However, such sign changes may lead to unexpected consequences, which will be explored further.

\vspace{0.1in}
{\noindent \bf Approximation of an Oscillatory Function.} 
We consider the oscillatory entire function \(f(x)=\sin(10\pi x)\) on \([-1,1]\). 
Since \(f\) is analytic, a degree-0 representation using Legendre polynomials converges 
exponentially once the approximating polynomial degree is sufficiently large to capture the 
oscillations. 
\begin{figure}[htbp]
    \centering
    \includegraphics[width=7.8cm, height = 4.2cm]{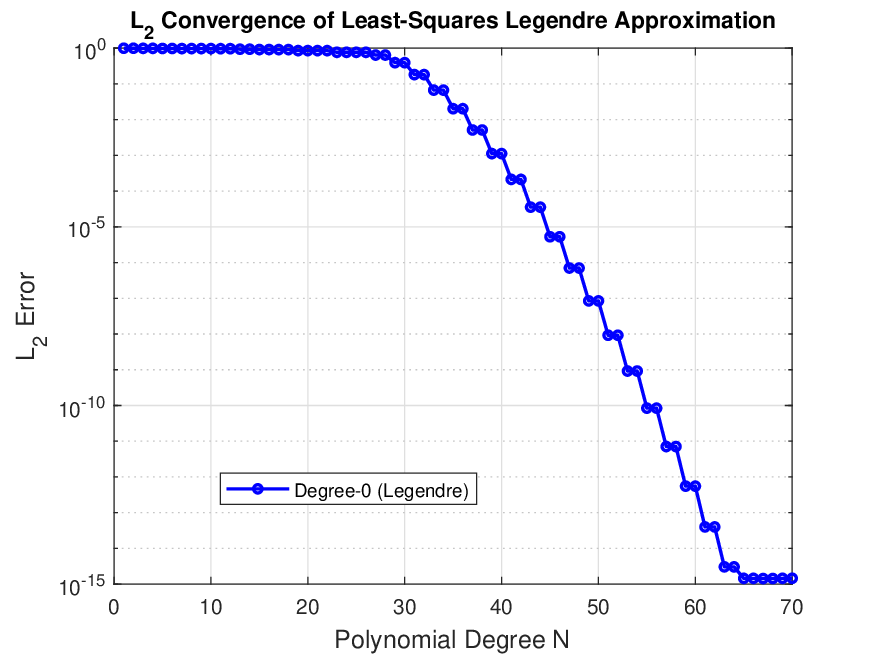}
    \caption{The \(L_2\) convergence of the Legendre approximation to $f(x)=\sin(10\pi x)$ as the polynomial degree $N$ increases. 
    After reaching a resolution threshold ($N\approx 30$), the error decays exponentially to machine precision.
    \label{fig:legendre-sin10-conv} }
\end{figure}
Figure~\ref{fig:legendre-sin10-conv}  illustrates the convergence of the \(L_2\) error
$$ \|f-f_N\|_{L_2}=\Big(\int_{-1}^1 |f(x)-f_N(x)|^2\,dx\Big)^{1/2}, $$
where $f_N(x)$ denotes the degree-0 approximation constructed using $N$ basis functions. 
Initially, the error decays very slowly when $N <30$.
In particular, when $N<4$, degree-0 approximations provide almost no useful 
information about $f(x)$. Beyond $N=30$, convergence transitions 
to exponential decay and ultimately reaches machine precision. This behavior is consistent with the theoretical convergence results for Legendre polynomial approximations. 

\begin{figure}[htbp]
    \centering
    \includegraphics[width=8.8cm,height=5.2cm]{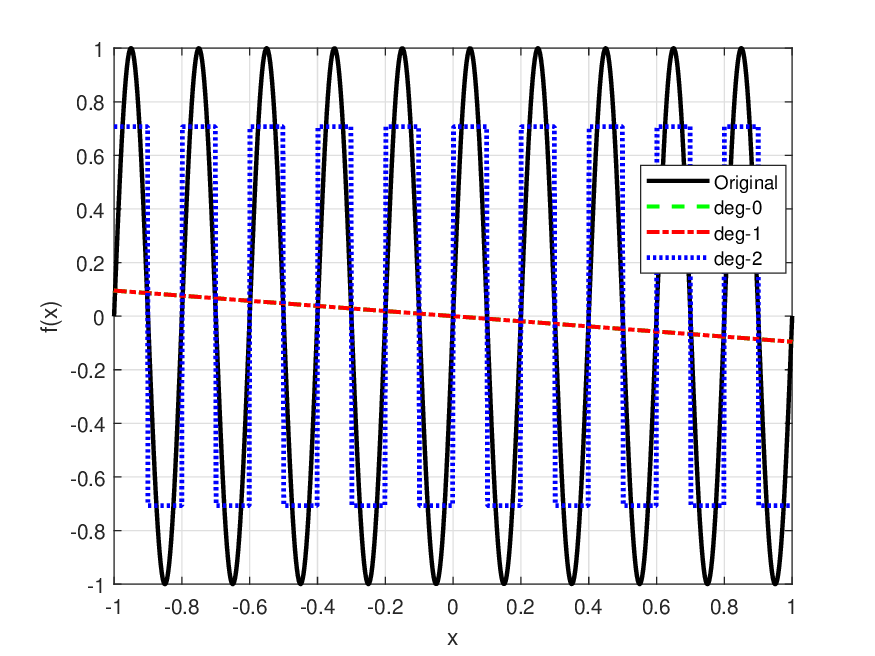}
    \caption{Point-wise approximations of \(f(x)=\sin(10\pi x)\) on \([-1,1]\) using two basis functions. Note that the Degree-0 representation (green dashed line) is identical to the degree-1 representation(red dash-dotted line).
    \label{fig:oci-comparison} }
\end{figure}

Fig.~\ref{fig:oci-comparison} compares the degree-0, degree-1, and degree-2 representations when only two basis functions are selected using a greedy algorithm to solve the least-squares optimization formulation. In this case, the degree-0 representation is a linear approximation, $f(x)=c_0+c_1 x$, with $c_0 \approx 0$ and $c_1 \approx -0.095$. The degree-1 approximation is either identical to the degree-0 approximation (when both basis functions are from $S_1$), or the least-squares solutions of $\| f(x) - c_0 -b_1 x f(x) \|$ 
(one basis function from $S_1$ and the other from $S_2$), or the solution of $\| f(x)  -b_1 x f(x) -b^2 x^2 f(x) \|$ (when both basis functions are from $S_2$). In our numerical experiment, the current greedy algorithm implementation produces the same approximation as the degree-0 representation.
The manifold corresponding to the optimal degree-2 representation is described by the quadratic equation
$f^2 - \frac{1}{2}=0$ corresponding to either $\{ c_0=\frac{1}{2}, c_1=0 \}$, 
or $\{ c_0=\frac{1}{2}, b_0=0 \}$, or $\{ c_0=\frac{1}{2}, a_1=0\}$. This representation provides a good approximation of the continuous wave using a discrete wave, and the periodicity of the original function is correctly captured by the index function $\zeta(x)$. 

\begin{figure}[htbp]
    \centering
    \includegraphics[width=7.2cm, height = 5.2cm]{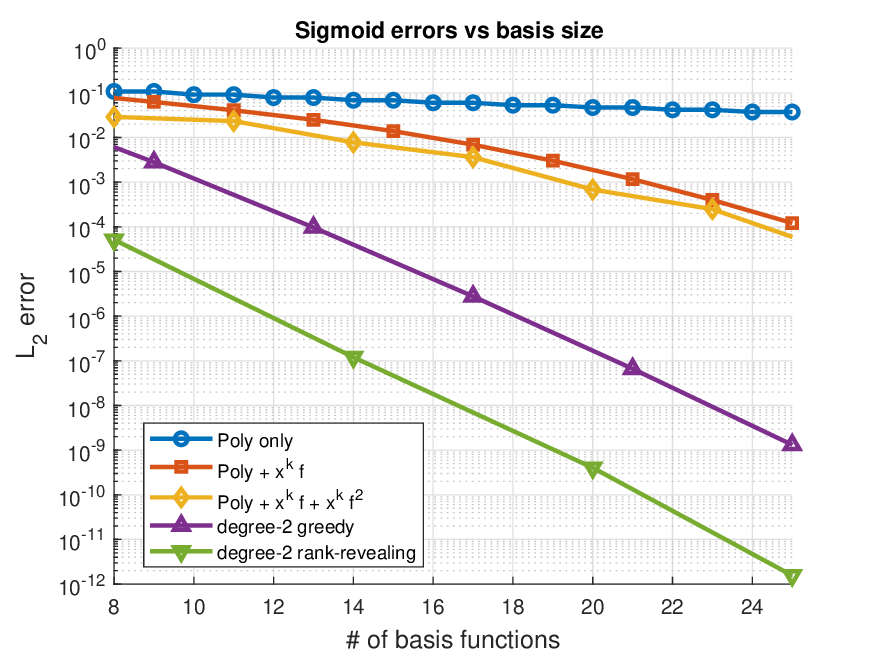}
    \caption{\(L_2\) error of the sigmoid approximation versus model size \(K\).
    Degree–2 adaptive rank-revealing basis selection technique outperforms all other approximation
    schemes.}
    \label{fig:sigmoid-conv}
\end{figure}

\vspace{0.1in}
{\noindent \bf Approximation of the Sigmoid Function with a Sharp Transition.} We consider the Sigmoid function
$
f(x)=\frac{1}{1+e^{-60x}}$, $ x\in[-1,1],
$
which exhibits a rapid transition near \(x=0\). 
Fig.~\ref{fig:sigmoid-conv} illustrates how different representations converge when the same number of basis functions is used. Although the function is analytic, the degree-0 representation (\textcolor{blue}{blue-$\circ$}) requires a large number of basis functions to resolve the sharp transition region. The non-adaptive degree-1 representation (\textcolor{red}{red-$\square$}, using equal number of basis functions from $S_1$ and $S_2$) converges more rapidly than the degree-0 representation due to its larger analyticity radius. The non-adaptive degree-2 representation (\textcolor{orange}{orange-$\lozenge$}, $N_0=N_1=N_2$) achieves slightly better convergence than the non-adaptive degree-1 representation. Degree-2 representations constructed via the greedy algorithm (\textcolor{violet}{violet-$\vartriangle$}, adding $s=1$ new basis at each step) exhibit exponential error decay. The best approximation is achieved using rank-revealing QR-filtered basis functions (\textcolor{green}{green-$\triangledown$}).

\vspace{0.1in}
{\noindent \bf Comment on the Stability of Degree-0, Degree-1, and Degree-2 Representations.} 
Thanks to advances in orthogonal linear algebra, the degree-0 representation can be constructed stably using the least-squares technique based on QR decomposition. However, constructing rational function approximation presents several challenges. One efficient and stable approach is the AAA algorithm \cite{nakatsukasa2018aaa}.
It is based on the rational barycentric representation, which provides a type $(m-1,m-1)$ rational approximation (degrees of the numerator and denominator polynomials) to the data set $\{ f(z_i) \}_{i=1}^m$. Although the AAA algorithm does not produce the optimal $L_2$ or $L_\infty$ rational approximation, it can still provide a useful representation for many practical applications when the sample points $\{ z_i \} _{i=1}^m$ are appropriately chosen. As a generalization of the degree-0 and degree-1 representations, the new degree-2 representation will also face several computational challenges, including efficient and accurate computation and numerical stability. We are particularly interested in determining the choices of $N_0$, $N_1$, and $N_2$, or their equivalent basis selections in the least-squares formulation. An adaptive strategy should be introduced to determine these choices. When the 
function is smooth, a degree-0 representation is preferable if it meets the requirements of the intended applications. However, for functions with discontinuities, degree-2 representations may offer advantages, as demonstrated in Sec.~\ref{sec:denoise}. In our numerical
experiments, we observed that for a piecewise-smooth function with a
single jump (see left panels in Fig.~\ref{fig:figJumps}), 
setting $N_0=N_1=N$ and $N_2=0$ yields a stable and accurate representation once $N$ is sufficiently large to fully resolve the manifold equation $f^2(x)-b(x) f(x) -c(x)=0$. A more detailed analysis of the properties of the degree-2 representation is currently underway.

\section{Data Denoising Using Degree-2 Representation}
\label{sec:denoise}
The separation of the smooth-coefficient manifold equation $a(x) f^2(x) - b(x) f(x) - c(x)=0$ and the index function
$\zeta(x)$ in Definition~\ref{def:deg2}, where $x \in [x_{\min}, x_{\max}]$ and 
$a(x)$, $b(x)$, and $c(x)$ are low-degree polynomials, 
enables more efficient processing of data with discontinuities, such as edge and corner singularities in computer images. 

Edge-preserving priors such as total variation (TV) \cite{rudin1992nonlinear} and sparse $L_1$ penalties \cite{donoho2006compressed} can be designed, however, they often lead to more challenging optimization problems. In this section, focusing on a single-variable function with jumps, we demonstrate how the degree-2 representation enables effective data denoising via least-squares-type optimization on the smooth-coefficient manifold equation $a(x) f^2(x) - b(x) f(x) - c(x)=0$. When noise is small, the index function can be determined directly 
by the distance between the data and the two branches. For large noise, statistical or network-based methods can be used to determine whether a point $x$ belongs to the first community 
(root 1) or the second community (root 2). As we focus on piecewise-smooth functions, we assume that the ground-truth image data lie exactly on a quadratic manifold
\begin{equation}\label{eq:deg2-manifold}
f^2(x)-b(x)f(x)-c(x)=0,
\end{equation}
with low-degree polynomials $b(x)=\sum_{n=0}^{N}b_n x^n$ and 
$c(x)=\sum_{n=0}^{N}c_n x^n$. A single jump is assumed at $x_0$; that is, the index function $\zeta(x)$ equals $-1$ for 
$x\in [x_{\min},x_0]$ and $1$ for $x \in (x_0,x_{\max}]$.
Although the degree-2 representation can be used for functions with multiple well-separated jumps, it may not be as effective as higher-degree representations. For instance, a degree-3 representation can readily handle piecewise-smooth functions with two jumps (3 roots). The ideas presented in this paper can be generalized to higher-degree representations in higher dimensions. These generalizations, currently under investigation, are discussed in
Sec.~\ref{sec:conclusions}.

In the numerical experiments in this section, we set the exact manifold equation 
as $(f-25)(f-255)=0$, where data are sampled at integer points on [0, 400], and the index function is defined as $-1$ when $0 \leq x \leq 140$ 
(root $f=25$) and $1$ when $140<x \leq 400$ (root $f=255$). Instead of an algebraic variety in the form $\hat{f}^2-b_0 \hat{f}-c_0=0$, we search for an approximate manifold equation $\hat{f}^2-(b_0+b_1 x) \hat{f}-(c_0+c_1 x)=0$ to simulate the case where more than enough terms in the expansions of $a(x)$, $b(x)$, and $c(x)$ are used to fully resolve the manifold.

Effective data denoising requires knowledge of the background noise. We consider four noise scenarios: (1) only that the collected data $\tilde f$ contains small noise (relative to the jump height), with no prior knowledge of other statistical properties; (2) white noise $\varepsilon(x)$ with mean $0$ added to the manifold equation, where the collected data $\tilde f$ satisfies 
$(\tilde f-25)(\tilde f-255)=\varepsilon$; (3) white noise $\varepsilon(x)$ with mean $0$ and known variance $\sigma^2$ added to the function, so that $\tilde f=f+\varepsilon$; and (4) a more general case where only certain statistical 
properties of the noise are known, such as its dependence on $f$ and interaction with different basis polynomials (e.g., statistical moments).

\vskip0.1in
\noindent\textbf{Case 1: Small but Otherwise Unknown Noise.} Given the observed data $\tilde f $ at each location, we first consider the case
where $\tilde f = f + \varepsilon$. We assume that $\| \varepsilon \|$ is small relative to the 
jump height $255-25=230$. However, no additional statistical properties of $\varepsilon$ are available; e.g., its mean, variance, and dependence on $f$ are unknown. Under these conditions, we estimate the parameters $(b_0,b_1,c_0,c_1)$ by solving the least-squares problem 
\begin{equation} 
\label{eq:case1optim}
  \min_{\{b_0,b_1,c_0,c_1\}} \| \tilde{f}^2-(b_0+b_1 x) \tilde{f}-(c_0+c_1 x) \|_2.
\end{equation}
The degree-2 representation is then given by the manifold equation 
$f^2-(b_0+b_1 x) f-(c_0+c_1 x)=0$, and the index function $\zeta(x)$ is defined by selecting the root closest to the sampled data $\tilde f(x)$ at each location $x$.

\begin{figure}[htbp]
\centering
\includegraphics[scale=0.5]{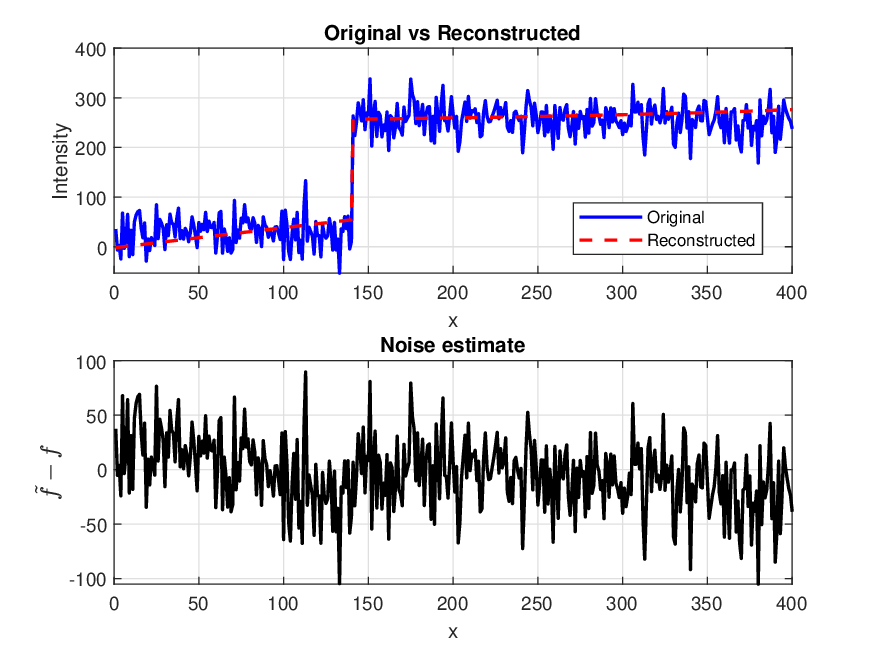}
\caption{Small error $\varepsilon$ is added to the function: $\tilde f=f+\varepsilon$ with 
  $\varepsilon=30\cdot\tt{randn}$. Top: original data $\tilde f$ (\textcolor{blue}{blue solid line}) and recovered $f$ (\textcolor{red}{red dashed line}).
  Bottom: estimated noise $\tilde{f}-f$.}
\label{fig:function-noise}
\end{figure}

For a numerical experiment, we set $\tilde f=f+\varepsilon$,
where $\varepsilon=30\cdot\tt{randn}$, and $\tt randn$ is the MATLAB function that generates random numbers from a normal distribution with mean 0 and variance 1.
The upper plot of Fig.~\ref{fig:function-noise} shows $\tilde{f}$ (\textcolor{blue}{blue solid line}) and the reconstructed $f$ (\textcolor{red}{red dashed line}). However, since the algorithm is unaware that the noise is white independent of $f$, the lower plot of Fig.~\ref{fig:function-noise} shows that the estimated noise $\tilde f - f$ is not white noise independent of $f$. The expected value $\mathbb{E}(\tilde f - f)$ is 
positive for $0 <x \leq 60$ and negative for $360 \leq x <400$. The estimated error remains within the ``small'' error regime.

\begin{figure}[htbp]
\centering
\includegraphics[scale=0.4]{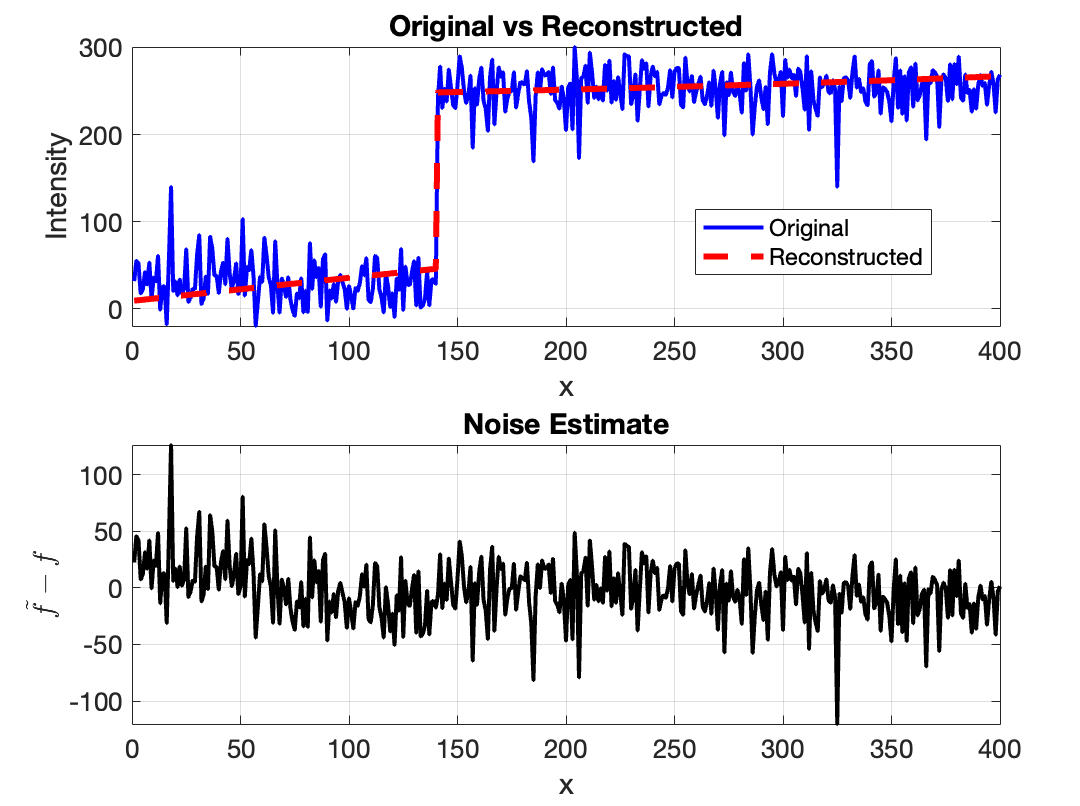}
\caption{White noise added to the manifold: $(\tilde f-25)(\tilde f-255)=\varepsilon$ with $\varepsilon=5000\cdot\tt{randn}$. Top: noisy data $\tilde f$ (\textcolor{blue}{blue solid line}) and recovered $f$ (\textcolor{red}{red dashed line}).
  Bottom: estimated noise $\tilde{f}-f$.}
\label{fig:manifold-noise}
\end{figure}

\vskip0.1in
\noindent\textbf{Case 2: White Noise with Zero Mean Added to Manifold Equation.} We assume white noise is added to the manifold equation
$(\tilde f-25)(\tilde f-255)=\varepsilon$, where $\varepsilon=5000\cdot\tt{randn}$. The noisy data $\tilde f$ is recovered using the quadratic formula, and the manifold equation is constructed by solving the least-squares problem: 
$\min_{\{b_0,b_1,c_0,c_1\}}\| \tilde{f}^2-(b_0+b_1 x) \tilde{f}-(c_0+c_1 x) \|_2$. The upper plot of Fig.~\ref{fig:manifold-noise} shows the noisy data $\tilde f$ (\textcolor{blue}{blue solid line}) alongside the reconstructed $f$ (\textcolor{red}{red dashed line}). Because the noise is relatively small, the index function $\zeta(x)$ is defined by selecting the root closest to the sampled data $\tilde f(x)$. The degree-2 representation effectively captures the original ground-truth function. 
The lower plot of Fig.~\ref{fig:manifold-noise} shows the estimated noise $\tilde f - f$. Although $\varepsilon$ is white noise, the estimated $\tilde f-f$ becomes colored noise, and its statistical properties depend on $f$. 

\vskip0.1in
\noindent\textbf{Case 3: White Noise $\varepsilon$ with Zero Mean and Known Variance $\sigma^2$ Added to 
the Ground-Truth Function $\tilde f = f + \varepsilon$.} In this case, we assume white noise is added to the ground-truth function: 
$\tilde{f} = f + \varepsilon$. The white noise has zero mean at each location, a known variance $\sigma^2$, and is independent of the function $f(x)$. 
Let $r(x)=f^2(x) - (b_0+b_1 x) f(x) - (c_0+c_1 x)$. Then, the least-squares solution of $\min_{\{b_0,b_1,c_0,c_1\}}\| f^2-(b_0+b_1 x) f-(c_0+c_1 x) \|_2$
satisfies
\begin{equation}
\langle r(x),1\rangle=\langle r(x),x\rangle=\langle r(x),f\rangle=\langle r(x),x f\rangle=0,\label{eq:ortho}
\end{equation}
where $\langle \cdot, \cdot \rangle$ denotes the standard inner product of two functions. This implies that $r(x)$ lies
in the subspace orthogonal to the basis functions $\{1, x, f(x), xf(x) \}$.
Let 
\[
S_0=\langle 1,1\rangle,\qquad S_x=\langle x,1\rangle,\qquad S_{x^2}=\langle x^2,1\rangle, 
\ \mbox{and} \ \ m_{x^p f^q}=\langle x^p, f^q\rangle,
\]
where $p=0,1,2$ and $q=1,2$. Then, enforcing the orthogonality condition in Eq.~(\ref{eq:ortho}), the coefficients $(b_0,b_1,c_0,c_1)$ are found by solving the 
linear system 
\begin{equation}\label{eq:lin-system}
\begin{bmatrix}
m_f & m_{x f} & S_0 & S_x \\
m_{x f} & m_{x^2 f} & S_x & S_{x^2} \\
m_{f^2} & m_{x f^2} & m_f & m_{x f} \\
m_{x f^2} & m_{x^2 f^2} & m_{x f} & m_{x^2 f}
\end{bmatrix}
\begin{bmatrix} b_0\\ b_1\\ c_0\\ c_1\end{bmatrix}
=
\begin{bmatrix}
m_{f^2}\\ m_{x f^2}\\ m_{f^3}\\ m_{x f^3}
\end{bmatrix}.
\end{equation}
Note that $S_0$, $S_x$, and $S_{x^2}$ can be computed easily. However, $m_{x^p f^q}$ cannot be computed directly since $f$ is unknown. Using the known properties of noise, $m_{x^p f^q}$ is derived from ``noisy'' moments 
\[
\tilde m_{x^p f^q}:=\langle x^p, \tilde f^q\rangle,~ p=0,1,2 \mbox{ and } q=1,2,
\]
where $\tilde{f}$ is the given noisy data.

\vspace{0.1in}
{\noindent \bf Denoising the Moments. }
Using the properties of the noise: $\mathbb{E}[\varepsilon]=0$, $\mathrm{Var}(\varepsilon)=\sigma^2$, and $\varepsilon(x)$ is independent of $f(x)$, we derive the expected relationship between the true moments $m_{x^p f^q}$ and the computed moments $\tilde m_{x^p f^q}$:
\begin{align*}
&\tilde m_f
=\langle 1,\tilde f\rangle
=\langle 1,f\rangle+\langle 1,\varepsilon\rangle
\stackrel{\mathbb{E}}{=} m_f,\\[2pt]
&\tilde m_{x f}
=\langle x,\tilde f\rangle
=\langle x,f\rangle+\langle x,\varepsilon\rangle
\stackrel{\mathbb{E}}{=} m_{x f},\\[2pt]
&\tilde m_{x^2 f}
=\langle x^2,\tilde f\rangle
=\langle x^2,f\rangle+\langle x^2,\varepsilon\rangle
\stackrel{\mathbb{E}}{=} m_{x^2 f},\\[4pt]
&\tilde m_{f^2}
=\langle 1,(f+\varepsilon)^2\rangle
=\langle 1,f^2\rangle+2\langle 1,f\varepsilon\rangle+\langle 1,\varepsilon^2\rangle
\stackrel{\mathbb{E}}{=} m_{f^2}+\sigma^2 S_0,\\[2pt]
&\tilde m_{x f^2}
=\langle x,(f+\varepsilon)^2\rangle
=\langle x,f^2\rangle+2\langle x,f\varepsilon\rangle+\langle x,\varepsilon^2\rangle
\stackrel{\mathbb{E}}{=} m_{x f^2}+\sigma^2 S_x,\\[2pt]
&\tilde m_{x^2 f^2}
=\langle x^2,(f+\varepsilon)^2\rangle
=\langle x^2,f^2\rangle+2\langle x^2,f\varepsilon\rangle+\langle x^2,\varepsilon^2\rangle
\stackrel{\mathbb{E}}{=} m_{x^2 f^2}+\sigma^2 S_{x^2},\\[4pt]
&\tilde m_{f^3}
=\langle 1,(f+\varepsilon)^3\rangle
=\langle 1,f^3\rangle+3\langle 1,f^2\varepsilon\rangle+3\langle 1,f\varepsilon^2\rangle+\langle 1,\varepsilon^3\rangle
\stackrel{\mathbb{E}}{=} m_{f^3}+3\sigma^2 m_f,\\[2pt]
&\tilde m_{x f^3}
=\langle x,(f+\varepsilon)^3\rangle
=\langle x,f^3\rangle+3\langle x,f^2\varepsilon\rangle+3\langle x,f\varepsilon^2\rangle+\langle x,\varepsilon^3\rangle
\stackrel{\mathbb{E}}{=} m_{x f^3}+3\sigma^2 m_{x f}.
\end{align*}
Using the known variance $\sigma^2$, we recover the true moments
\begin{align*}
m_f          &= \tilde{m}_f, & 
m_{xf}       &=\tilde{m}_{xf}, & 
m_{x^2 f}    &= \tilde{m}_{x^2 f}, \\
m_{f^2}      &= \tilde m_{f^2} - \sigma^2 S_0, &
m_{x f^2}    &= \tilde m_{x f^2} - \sigma^2 S_x, &
m_{x^2 f^2}  &= \tilde m_{x^2 f^2} - \sigma^2 S_{x^2},\\
m_{f^3}      &= \tilde m_{f^3} - 3\sigma^2 \tilde m_f, &
m_{x f^3}    &= \tilde m_{x f^3} - 3\sigma^2 \tilde m_{x f}.
\end{align*}
Finally, solving the linear system in Eq.~(\ref{eq:lin-system}) yields the denoised manifold equation. 

\vspace{0.1in}
{\noindent \bf Denoising the Index Function $\zeta(x)$. }
When the noise variance is small, the index function $\zeta(x)$ can be 
defined by selecting the root closest to the sampled data $\tilde f(x)$.
However, if the noise at a particular location $x$ becomes large relative to the jump height, the observed data $\tilde f$ may appear closer to a different root. In such cases, this simple scheme may return an incorrect index function value $\zeta(x)$. To denoise $\zeta(x)$, 
we introduce a k-nearest neighbor (k-NN)-based collective ``voting'' strategy to determine the correct index value. For each point $x$, we first identify its k-NN, where $k$ is a pre-determined parameter that depends on the variance $\sigma^2$. We set
$\zeta(x)$ to the root selected by the majority of $x$ and its k-NNs. We repeat the voting process until $\zeta(x)$ stabilizes.

\begin{figure}[htbp]
\centering
\includegraphics[scale=0.4]{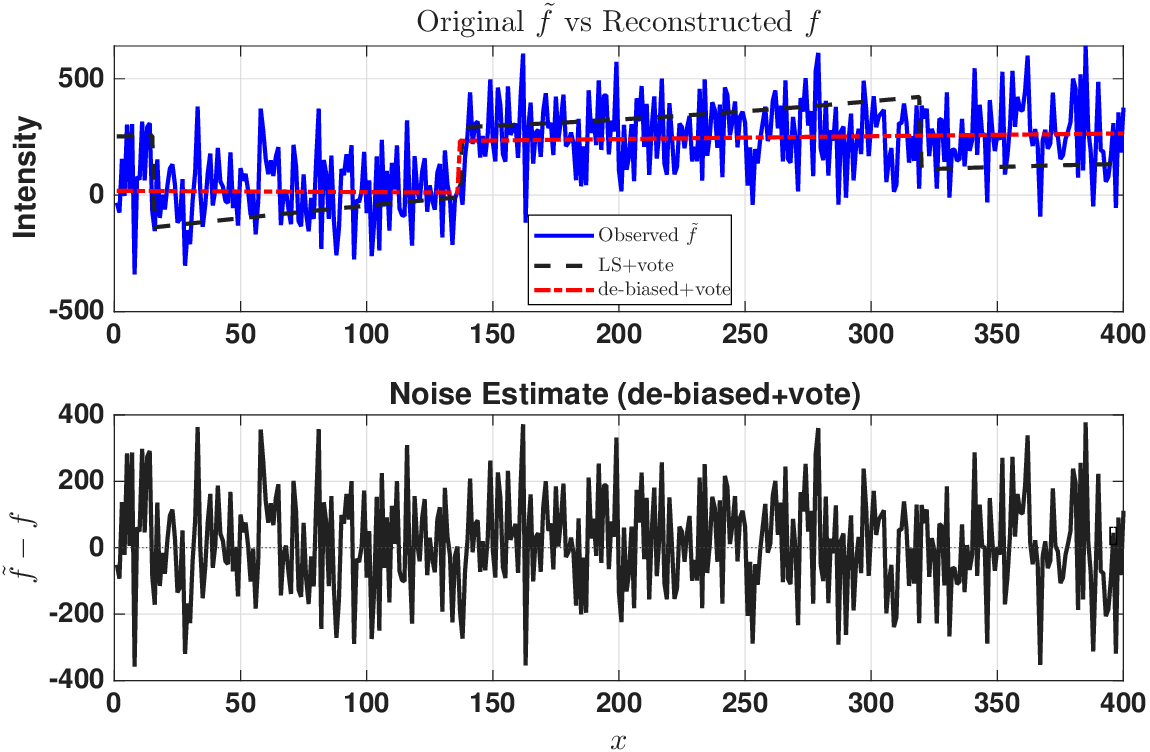}
\caption{High-variance function noise with known variance: $\tilde f=f+\varepsilon$, $\varepsilon=150\cdot\tt{randn}$, $\sigma^2=150^2$ used for de-biasing. Top: observed $\tilde f$ (\textcolor{blue}{blue solid line}), LS+vote (no correction)(black dashed line), and de-biased+vote (known-variance correction)(\textcolor{red}{red dashed-dotted}). Bottom: estimated noise.}
\label{fig:function-noise-highvar}
\end{figure}

The upper plot of Fig.~\ref{fig:function-noise-highvar} shows three graphs: (1) the observed noisy data $\tilde f$ (\textcolor{blue}{blue solid curve}); (2) the denoised degree-2 function $f$ obtained by deriving the manifold equation using the least-squares
formulation in Eq.~(\ref{eq:case1optim}) in Case 1 and applying the voting strategy (\textcolor{black}{black dashed curve}, LS+vote); and (3) a further improved estimate of $f$ obtained by solving Eq.~(\ref{eq:lin-system}) with revised moments, combined with the voting strategy (\textcolor{red}{red dash-dotted curve}, de-biased+vote). 
Using the statistical properties of the noise and de-biased moments, the resulting degree-2 representation (\textcolor{red}{red dash-dotted line}) effectively recovers the ground-truth function $f$. 
The lower plot of Fig.~\ref{fig:function-noise-highvar} shows the estimated noise $\tilde{f}-f$, which is approximately white and passes the independence statistical test ($\tilde{f}-f$ is independent of $x$).

\vskip0.1in
\noindent\textbf{Case 4: More General Noise.} Finally, we consider the case where only limited statistical properties of the noise are known, such as
its dependence on the function $f(x)$. We propose the following iterative approach to improve the
accuracy of the denoised function. 

\vspace{0.1in}
\begin{enumerate}
\item \textbf{Initialization.} Using the observed noisy data, determine the initial degree-2 representation functions $a_0(x)$, $b_0(x)$, $c_0(x)$ and the index function $\zeta_0(x)$ by applying the method from Case 1 or Case 2, or by using the denoising schemes from Case 3 with an initial variance estimate $\hat\sigma^2_0$. 

\item \textbf{For $i=0,1,2,\dots$}
\begin{enumerate}
\item \textbf{Compute the estimated noise $\tilde{\varepsilon}_i=\tilde f - f_i$:}  

Using the estimated degree-2 coefficient functions $a_i(x)$, $b_i(x)$, $c_i(x)$, and the index 
function $\zeta_i(x)$, construct $f_i$ and compute the 
estimated noise $\tilde{\varepsilon}_i =\tilde f -f_i$.
\item \textbf{Denoise $\tilde{\varepsilon}_i$ and obtain improved data $\bar{f} = \tilde{f}-\bar{\varepsilon}$:} 

We assume the estimated noise is the sum of the real noise $\varepsilon$ and a smooth function represented by a degree-0 expansion 
$$\tilde{\varepsilon}_i(x) = \varepsilon(x) + \sum_{n=0}^{K-1} c_n L_n(x),$$
where $L_n$ denotes the Legendre polynomial of degree $n$, $c_n$ are unknown coefficients,
and $K$ represents the number of conditions known a priori about the real noise $\varepsilon$, e.g., $\mathbb{E} \left[ \varepsilon(x)\right]=0$, $\langle x,\varepsilon\rangle =0$, etc.
Determine $c_n$ using these conditions. Then define the improved noise $\bar{\varepsilon}(x)$ as
$$\bar{\varepsilon}(x) =\tilde{\varepsilon}_i(x) - \sum_{n=0}^{K-1} c_n L_n(x).$$ 
This procedure can be interpreted as projecting the estimated noise onto a subspace defined by the available statistical
properties of the noise. Subtract the improved noise $\bar{\varepsilon}$ from the observed data $\tilde f$ and denote the improved denoised data as $\bar f = \tilde{f} - \bar{\varepsilon}$. 

\item \textbf{Compute de-biased  $a_{i+1}(x)$, $b_{i+1}(x)$, $c_{i+1}(x)$, and $\zeta_{i+1}(x)$:}  

Using the improved denoised data $\bar f$ from (b), formulate and solve the least-squares problem 
to obtain the updated, de-biased degree-2 coefficient functions
 $a_{i+1}(x)$, $b_{i+1}(x)$, $c_{i+1}(x)$, and the index function $\zeta_{i+1}(x)$. 
 
 \item \textbf{Check convergence:}
 
If the new coefficient functions $a_{i+1}(x)$, $b_{i+1}(x)$, $c_{i+1}(x)$, and the index 
function $\zeta_{i+1}(x)$ have converged, return $a(x) = a_{i+1}(x)$, $b(x)=b_{i+1}(x)$, $c(x)=c_{i+1}(x)$, the index function $\zeta(x)=\zeta_{i+1}(x)$, and the reconstructed $f(x)$. Otherwise, repeat Step (a).
\end{enumerate}
\end{enumerate}
\begin{figure}[htbp]
\centering
\includegraphics[scale=0.4]{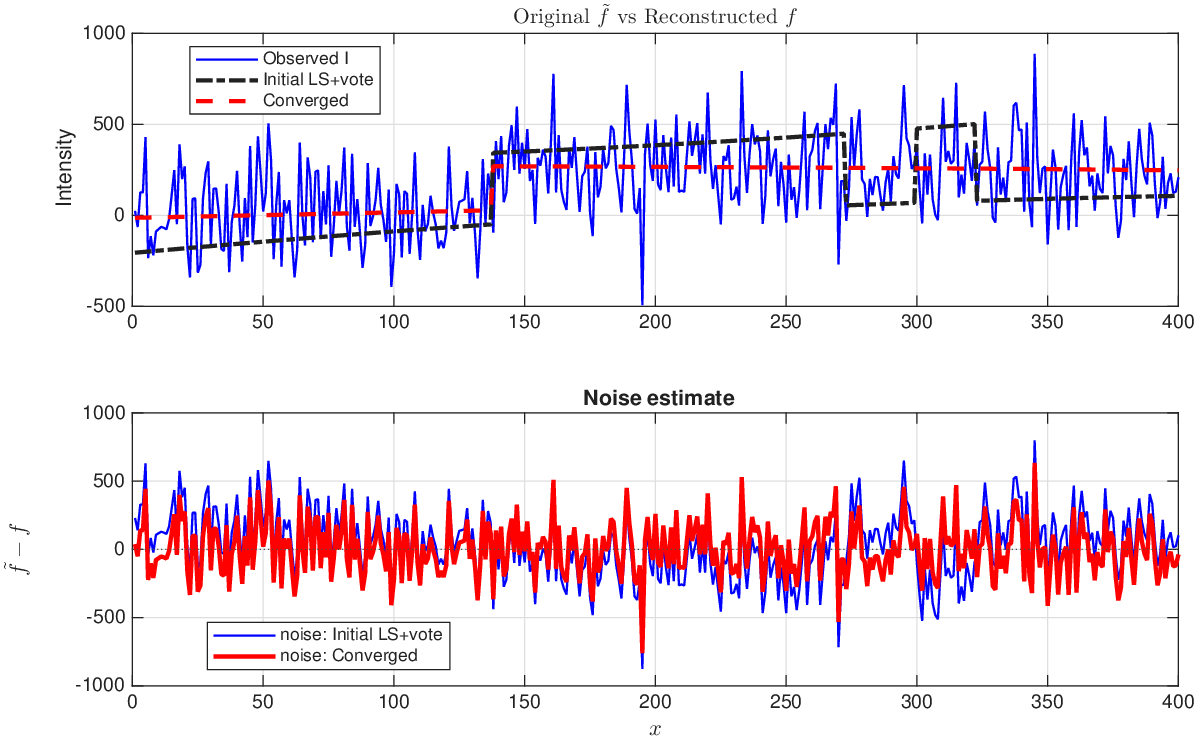}
\caption{High-variance function noise with unknown variance: $\tilde f=f+\varepsilon$, $\varepsilon=200\cdot\tt{randn}$. 
Top panel: observed $\tilde f$ (\textcolor{blue}{blue solid line}), initial reconstruction $f_0$ (black dash-dotted), and converged solution $f_{final}$(\textcolor{red}{red dashed line}) from the iterative de-biasing and voting technique. Bottom panel: estimated noise of the initial reconstruction (\textcolor{blue}{blue solid line}) and final reconstruction (\textcolor{red}{red solid line}).}
\label{fig:generalIterative}
\end{figure}

\vspace{0.1in}
In Fig.~\ref{fig:generalIterative}, we present the numerical results produced by the
iterative approach. The observed data (\textcolor{blue}{blue solid line}) $\tilde f=f+\varepsilon$, where 
$\varepsilon=200\cdot\tt{randn}$, was generated by adding white noise with a variance of $200^2$
to the ground-truth function $f$.
The initial reconstruction, denoted by $f$ (assuming small noise and solving 
the manifold equation via least-squares, followed by index/label voting), is shown as the 
black dash-dotted curve labeled ``LS+vote''. The noise was constrained to satisfy eight moment conditions, namely $\langle 1,\varepsilon\rangle =0$,
$\langle x,\varepsilon\rangle =0$, $\langle x^2,\varepsilon\rangle =0$,
$\langle f,\varepsilon\rangle =0$, $\langle x f,\varepsilon\rangle =0$,
$\langle x^2 f,\varepsilon\rangle =0$, $\langle f^2,\varepsilon\rangle =0$,
and $\langle x f^2,\varepsilon\rangle =0$. 
The iterative correction scheme converged, and the final reconstruction of $f$
is shown as the \textcolor{red}{red dashed line}, labeled ``Converged''.  
As shown in the top panel of Fig.~\ref{fig:generalIterative}, the converged solution closely matches the ground-truth function $f$. The bottom panel of Fig.~\ref{fig:generalIterative} shows the estimated noise both at initialization and after convergence. For the LS+vote baseline (\textcolor{blue}{blue solid line}) at initialization, the estimated noise exhibits significant bias: for example, in the first 100 indices, the noise remains predominantly positive, whereas around indices 200--260, it becomes more negative, revealing a low-frequency 
bias dependent on $x$. In contrast, after applying the iterative de-biasing corrections (\textcolor{red}{red solid line}), the estimated noise fluctuates symmetrically around zero across the entire domain. This behavior is consistent with the enforced moment constraints and the assumption that the true noise $\varepsilon(x)$ is 
independent of both $x$ and $f(x)$. These results demonstrate that the proposed iterative correction effectively eliminates smooth and low-frequency bias, thereby improving the accuracy of the recovered function $f$.
The convergence properties of the algorithm are still under investigation, and additional assumptions may be required to guarantee convergence.

\section{Conclusions and Future Work}
\label{sec:conclusions}
The main contribution of this paper is introducing a quadratic-formula-based
degree-2 nonlinear approximation scheme. The new presentation extends the classical degree-0 
and degree-1 representations by decoupling the smooth-coefficient manifold equation from the index function, which indicates the root or branch to select. The degree-2 representation offers a
more accurate approximation for functions with discontinuities or sharp transitions and is applicable to the analysis of multi-valued functions. Numerical algorithms were developed to construct degree-2 representations, and numerical experiments demonstrate their effectiveness and compare the proposed nonlinear degree-2 representation with existing degree-0 and degree-1 methods. As an illustrative application, the degree-2 representation is used for data denoising by independently processing the smooth-coefficient manifold equation and the index function. Numerical results for various noise types validate the effectiveness of the approach. 
The degree-2 representation belongs to an underexplored class of nonlinear approximation techniques, and its 
analytical and numerical properties warrant further investigation. In the 
following, we briefly outline several topics currently under investigation. 

Firstly, we obtained promising results with the current basis selection 
strategies described in Sec.~\ref{sec:deg2-construction}. We are exploring ways to enhance the construction process 
by integrating concepts from rank-revealing QR algorithms and other manifold learning 
techniques. For example, rather than identifying the most ``significant'' basis vectors, 
the presence of an algebraic relation suggests searching for a small set of basis 
functions that are nearly linearly dependent. However, the relation equation 
with the fewest basis terms may face numerical stability issues, since computing the roots of a general polynomial can be ill-conditioned. 
Therefore, selecting basis functions from the 
frame set in an appropriate order is essential for improved accuracy, efficiency, and stability. 

Secondly, the degree-2 representation can be generalized to higher degrees by fitting 
an algebraic relation between the data $x$ and $f(x)$ using
\begin{equation}\label{eq:degd-model}
\sum_{j=0}^{d} a_j(x)\,(f(x))^j \;\approx\; 0,
\qquad
a_j(x)=\sum_{n=0}^{N_j} \alpha_{j,n}\,x^n,
\end{equation}
where a simple normalization strategy is to set the highest-degree coefficient function to $a_d(0) =1$ or $a_d(x) \equiv 1$ (or $N_d=0$).
Once the degree-$d$ polynomial equation representation, or algebraic variety, $P(x,f(x))=0$
is constructed for an input variable $x$, we evaluate the corresponding coefficient 
functions $\{ a_j(x) \}_{j=0}^d$. Then, the function value $f(x)$ is obtained by solving
\begin{equation}\label{eq:degd-root}
P(x,z):=\sum_{j=0}^{d} a_j(x)\,z^{j}=0
\end{equation}
and using the index function to select the correct root or branch when $f$ is single-valued. Note that computing the roots of a polynomial is a well-studied problem \cite{boyd2002computing,zhang2024finding}. 
Both the greedy algorithm and rank-revealing basis selection techniques can be
extended to higher-degree representations. We are currently testing 
a two–axis greedy strategy to determine whether to increase the power/degree 
of $f$ or raise the coefficient degree for 
one of the existing power terms (by adding the next degree $x$–polynomial multiplied by an existing 
$f^{j}$ term). At each iteration, we compute the residual reduction for 
each candidate and select the one that achieves the greatest reduction.

Finally, this paper focuses on denoising single-variable functions. Extending the degree-$d$ representation to multivariable functions is straightforward. For example, a 2D image function can be represented as $\sum_{k=0}^K c_k(x,y) f^k(x,y) =0.$ 
The coefficients $c_k(x,y)$ are polynomials of the spatial variables $x$ and $y$. 
For video data, the representation extends to $\sum_{k=0}^K c_k(x,y,t) f^k(x,y,t) =0,$ which incorporates the temporal variable $t$ into the coefficient functions. We are currently exploring higher-dimensional degree-$d$ representations and their applications in image and video processing.

\section*{Acknowledgments}
We gratefully acknowledge the generous support from NSF grants DMS2152289 (J. Huang and Z. He) and 
DMS2152070 (Y. Wu). Part of this work was carried out during a Collaborate@ICERM group event, supported 
by the Institute for Computational and Experimental Research in Mathematics (ICERM) at Brown University. 
The authors thank ICERM for its hospitality and for fostering a stimulating,  collaborative research environment during
the event.

\bibliographystyle{siamplain}
\bibliography{references}
\end{document}